\documentclass[a4paper,12pt]{article}

\usepackage[latin1]{inputenc}
\usepackage[english]{babel}

\usepackage{amssymb}
\usepackage{amsmath}
\usepackage{amsthm}
\usepackage{amscd}
\usepackage{mathtools}
\usepackage{enumitem}

\usepackage{textcomp}
\usepackage{layout}
\usepackage{color}

\newcommand{\Poi}{\mathbf P} 

\newcommand{\Poy}{\mathsf S\!} 

\newcommand{\Pp}{\mathsf P} 

\newcommand{\B}{\mathcal{B}} 
\newcommand{\Bbd}{\mathcal{B}_{0}} 

\newcommand{\M}{\mathcal M}
\newcommand{\MX}{\mathcal M(X)}

\newcommand{\Mpm}{\mathcal M^{\cdot\cdot}}
\newcommand{\MpmX}{\mathcal M^{\cdot\cdot}(X)}

\newcommand{\N}{\mathbb N}

\newcommand{\R}{\mathbb R}

\renewcommand{\d}{\mathrm{d}}
\newcommand{\supp}{\operatorname{supp}}

\renewcommand{\exp}{\operatorname{exp}}

\newcommand{\vi}{\varphi}
\newcommand{\vo}{\varrho}

\theoremstyle{definition}
\newtheorem{defdefinition}{Definition}[section]

\theoremstyle{plain}
\newtheorem{defsatz}[defdefinition]{Theorem}
\newtheorem{defsatzdef}[defdefinition]{Theorem and Definition}
\newtheorem{defprop}[defdefinition]{Proposition}
\newtheorem{deflemma}[defdefinition]{Lemma}
\newtheorem{deffolgerung}[defdefinition]{Corollary}
\newtheorem{defbeispiel}[defdefinition]{Example}

\theoremstyle{remark}
\newtheorem{defbemerkung}[defdefinition]{Remark}

\newcommand{\defn}[2]{\begin{defdefinition}[#1]#2\end{defdefinition}}

\newcommand{\satz}[1]{\begin{defsatz}#1\end{defsatz}}
\newcommand{\satzn}[2]{\begin{defsatz}[#1]#2\end{defsatz}}

\newcommand{\prop}[1]{\begin{defprop}#1\end{defprop}}

\newcommand{\bem}[1]{\begin{defbemerkung}#1\end{defbemerkung}}
\newcommand{\lemma}[1]{\begin{deflemma}#1\end{deflemma}}

\newcommand{\korollar}[1]{\begin{deffolgerung}#1\end{deffolgerung}}

\newcommand{\beispiel}[1]{\begin{defbeispiel}#1\end{defbeispiel}}

\theoremstyle{plain}

\newcommand{\eM}{{\cal M}^{\cdot \cdot}}

\newcommand{\J}{{\cal J}}

\newcommand{\K}{{\cal K}}
\newcommand{\F}{{\cal F}}
\newcommand{\A}{{\cal A}}

\newcommand{\D}{{\cal D}}

\newcommand{\pe}{{\cal P}}

\title{The logical postulates of \\B\"oge, Carnap and Johnson\\in the context of Papangelou processes}

\author{Mathias Rafler\footnote{TU M\"unchen, Zentrum Mathematik M5, Boltzmannstr. 3, D-85747 Garching bei M\"unchen, Germany, email {\sf rafler@ma.tum.de}};
Hans Zessin\footnote{Fakult\"at f\"ur Mathematik der Universit\"at Bielefeld, Postfach 10 01 31,
\newline D-33601 Bielefeld, Germany, e-mail: {\sf zessin@math.uni-bielefeld.de} }
}

\begin{document}

\maketitle

\begin{abstract}
We adapt Johnson's sufficiency postulate, Carnap's prediction invariance postulate and B\"oge's learn-merge invariance to the context of Papangelou processes and discuss equivalence of their generalizations, in particular their weak and strong generalizations. This discussion identifies a condition which occurs in the construction of Papangelou processes. In particular, we show that these generalizations characterize classes of Poisson and P\'olya point processes.
\end{abstract}

\section{Introduction}

To construct locally finite, but possibly infinite random point configurations in abstract spaces, one might wih to build such a configuration point by point. In such a case, one needs conditional rules which state how to place a point given a point configuration. Such ``conditional intensities'' roughly give the intensity of realizing a point conditioned on having observed a particular configuration. We are interested in the structure of these objects and thus the point process itself: Postulating simple properties, how do this intensity and the point process look like? We give characterizations of several classes of point processes such as the class of Gibbs point processes and the class of P\'olya and Poisson point processes.

The point processes of our interest are point processes which admit a partial integration formula. They were firstly studied by Papangelou~\cite{P,P1} and are called Papangleou processes; the kernel that appears is the Papangelou kernel. This kernel allows to infer from an observed point configuration to an intensity for a new point to realize.

There are two main motivations to study the structure of Papangelou kernels: Firstly, the works of Johnson~\cite{J}, Carnap~\cite{CS58} and B\"oge and M\"ocks~\cite{BM86} on inductive reasoning and Bayesian learning models, respectively; secondly, a basic assumption used in the construction of Papangelou processes from Papangelou kernels in Zessin~\cite{hZ09}. We adapt these postulates to Papangelou kernels and identify the structure of these kernels and thus the corresponding Papangelou process.

The aim of Johnson and Carnap was to give a logical fundament to inductive reasoning, which means ``To infer from a given premise to an hypothesis, which logically goes beyond this premise, and thus such a reasoning can only be true with a certain probability.'' For a presentation of these ideas seeEssler~\cite{Ess}, Humburg~\cite{H}. Mathematically, starting from an exchangeable sequence of random variables, the aim is to infer from the observation of say the first $n$ variables on the conditional law of the next outcome. This law is in Carnap's language the inductive probability given the observations. By exchangeability, this probability should only depend on the counts of each possible outcome. Carnap suggests a prediction irrelevance postulate NA 14 in Carnap and Stegm\"uller~\cite{CS58}: The probability of some outcome of the $n+1$ observation should be invariant under altering the observations of any \emph{other} outcome. Johnson assumed that this probability shall only depend on the number of observations of that category, according to Good this is called sufficiency postulate~\cite{G}. A rigorous treatment of Johnson's reasoning is given by Zabell~\cite{sZ82}. Finally, B\"oge et al. postulate that it should make no difference in inference, whether one combines (merges) some of the observation categories before or after sampling. This principle they called learn-merge invariance and they showed that these concepts are equivalent in a certain framework. 

Their construction guarantees that, given an a priori probability as an initial condition, there exists a joint probability space for all observations and the postulates characterize the conditional laws given previous observations. Roughly, these are the laws of urn models with or without replacement or P\'olya's urn.

Basically, Johnson's and Carnap's postulates are measurability conditions: The intensity in some bounded set shall only depend on the number of points observed in that set so far. With the aid of a strong version of the sufficiency postulate, we characterize Poisson and P\'olya processes. The nature of the learn-merge invariance is, that, when passing from a micro- to a macrolevel, the transformation of the conditional probabilities is consistent. In our situation, mapping a Papangelou process with a given kernel yields a Papangelou process for the mapped kernel. We show that the consistency for a sufficiently large class of transformations is equivalent to the sufficiency postulate.

The big surprise is the following: The mentioned point processes represent the ideal gases of quantum statistical mechanics, Maxwell-Boltzmann, Bose-Einstein or Fermi-Dirac states. They are the only ones which satisfy these universal invariance properties. For details we refer to Bach and Zessin~\cite{BZ}.

During the construction of Papangelou processes by means of their kernels in Zessin~\cite{hZ09}, an absolute continuity condition relating the intensities for some configuration with the intensity for that configuration with an additional point was used. A special case is the basis of our argumentation and it turns out that the basic structure of these kernels is quite simple: A fixed measure plus some reinforcements depending on the configuration only in a local way. For the general condition and an additional assumption we identify the class of Gibbs point proceses.

In the following section, we give a short overview over Papangleou processes and the interpretation of the Papangelou kernel. In Section~\ref{sect:char}, we resolve the relation between the logical postulates and Papangelou kernels and characterize the class of P\'olya and Poisson processes. The follwoing section is devoted to the characterization of Gibbs point processes.

\section{Point Processes and Papangelou Processes}

\subsection{Point Processes}

Firstly we fix some notations, see~\cite{Ka,MKM78,LMW88}. Let $(X,d)$ be a complete, separable metric space (c.s.m.s.) with its Borel $\sigma$-field $\B(X)$ and the ring $\Bbd(X)$ of bounded Borel sets in $X$. For a measure $\mu$ on $X$ and a Borel set $B$ the evaluation mapping is denoted by $\zeta_B(\mu)=\mu(B)$. $\MX$ is the collection of Radon measures, i.e. of measures which are finite on $\Bbd(X)$. $\MX$ is again a c.s.m.s. when equipped with the topology generated by the collection of mappings $\{\zeta_B\}_{B\in\Bbd(X)}$. 

Denoting by $\zeta_f$ the mapping $\mu\mapsto\mu(f):=\int f\,\d\mu$, then the convergence of a sequence of measures $(\mu_k)_k$ in $\MX$ is equivalent to the convergence of $(\zeta_f(\mu_k))_k$ for each bounded, continuous function $f:X\to\R$ with bounded support. 

Let $\MpmX\subset\MX$ be the set of point measures, i.e. the set of measures $\mu\in\MX$ such that $\mu(B)\in\N\cup\{+\infty\}$ for all $B\in\B(X)$. $\MpmX$ is a measurable subset of $\MX$. A probability measure $\Pp$ on $\MX$ is called \emph{random measure}, and if $\Pp\bigl(\MpmX\bigr)=1$, $\Pp$ is called \emph{point process}. We also write $\Pp(\vi)$ for the integral $\int\vi\,\d\Pp$. Finally, let $\Mpm_n(X)$ be the set of point measures of total mass $n$ for each $n\in\N$.

By $F(X)$ denote the set of positive, measurable functions on $X$ and by $F_c(X)$ 
continuous functions with bounded support.

\subsection{Papangelou Processes}

The set $\{(x,\mu)\in X\times\MpmX: \mu(\{x\})>0\}$ is a measurable subset of the product space, and the Campbell measure $C_\Pp$ of a point process $\Pp$ is supported by this set and defined as
\begin{equation*}
  C_\Pp(h)=\iint h(x,\mu)\,\mu(\d x)\Pp(\d\mu),\qquad h\in F\bigl(X\times\MpmX\bigr).
\end{equation*}
If $\Pp$ is of first order, then $C_\Pp$ is a $\sigma$-finite measure. Moreover, $\Pp$ is uniquely determined by $C_\Pp$. If the function $h$ does not depend on the second component, then we obtain the intensity measure of $\Pp$. If this measure is Radon, $\Pp$ is said to be {\it of first order}.

A mapping $\pi:\MpmX\times\B(X)\to\R_+$ is a \emph{Radon kernel} from $\MpmX$ to $X$, if for each $\mu\in\MpmX$, $\pi(\mu,\,\cdot\,)\in\MX$ is a Radon measure on $X$, and moreover, the mapping $\mu\mapsto\pi(\mu,B)$ is measurable for each $B\in\B(X)$.

\defn{Papangelou Process}{
A \emph{Papangelou Process specified by the Radon kernel $\pi$} is a point process which satisfies the integration-by-parts formula
\begin{equation*}
   C_\Pp(h) = \iint_{X\times\MpmX} h(x,\mu+\delta_x)\,\,\pi(\mu,\d x)\Pp(\d\mu),\qquad h\in F .
\end{equation*}
In this case, $\pi$ is called \emph{Papangelou kernel}. Let $\wp(\pi)$ denote the collection of all Papangelou processes with Papangelou kernel $\pi$.
}

A Papangelou process determines uniquely the kernel $\pi$ up to some null set. The {\it probabilistic interpretation of a Papangelou process} is intimately related to the notion of a Gibbs state in classical statistical mechanics: Define the iterated kernels $\pi^{(k)}$ given recursively by
\begin{equation*}
  \pi^{(k)}(\mu;\d x_1,\ldots,\d x_k) = \pi(\mu+\delta_{x_1}+\cdots+\delta_{x_{k-1}},\d x_k)\pi^{(k-1)}(\mu;\d x_1,\ldots,\d x_{k-1}),
\end{equation*}
$\pi^{(1)}=\pi$ and $\pi^{(0)}=1$, and moreover
\begin{equation*}
  \pi_B(\mu,\vi) = \sum_{m=0}^\infty \frac{1}{m!} \int_{B^m} \vi(\delta_{x_1}+\cdots +\delta_{x_m})\,\,\pi^{(m)}(\mu_{B^c};\d x_1,\ldots,\d x_m).
\end{equation*}
If $\Pp$ is a Papangelou process of first order specified by $\pi$, then, see e.g.~\cite{MWM79,NZ79},
\begin{equation*}
  0 < \pi_B(\,\cdot\,,\eM) < \infty\qquad \Pp\text{-a.s}
\end{equation*}
 for each $B\in \B_0$. Moreover, if $\chi_B(\mu)$ is the restriction of a configuration $\mu$ to $B$, then the conditional distribution of $\chi_B$ given the environment outside $B$ is the normalization of the kernel $\pi_B$, i.e.
\begin{equation*}
  \Pp\bigl( \vi\circ\chi_B|\chi_{B^c}\bigr) = \frac{\pi_B(\,\cdot\,,\vi)}{\pi_B(\,\cdot\,,\Mpm)}\qquad \Pp\text{-a.s},
\end{equation*}
where $\vi\in F$. Thus, Papangelou processes are locally specified by the kernel $\pi$, given the environment.

A kernel $\pi$ is \emph{simple} if $\pi(\mu,\supp\mu)=0$ for all $\mu\in\MpmX$. Note that a Papangelou process $\Pp$ is simple if and only if $\pi$ is simple.

\beispiel{ \label{ex:pap:kerne}
\begin{enumerate}
  \item If $\pi(\mu,\,\cdot\,)=\vo$, then $\wp(\pi)$ consists of the Poisson process $\Poi_\vo$ with intensity measure $\vo$. This is known as Mecke's characterization of the Poisson process. Observe that $\pi$ is simple if and only if $\vo$ is diffuse.
  \item If $\pi(\mu,\,\cdot\,)=z(\vo+\mu)$ for some $z\in(0,1)$, then again $\wp(\pi)$ contains a single element, the P\'olya sum process $\Poy_{z,\vo}$. Clearly, $\pi$ is not simple.
  \item If $\pi(\mu,\,\cdot\,)=z(\vo-\mu)$, where $z>0$ and $\mu\leq\vo\in\MpmX$ and $\pi=0$ otherwise, then again $\wp(\pi)$ is a singleton and consists only of the P\'olya difference process. Moreover, $\pi$ is simple if and only if $\vo$ is a simple measure.
\end{enumerate}
}
For the second process see Zessin~\cite{hZ09}, for the third one~\cite{NZ11}. In general, $\wp(\pi)$ may also contain a continuum of Papangelou processes or even be empty. A necessary condition for a Radon kernel to be a Papangelou kernel is the cocycle condition given in $(\A_1)$ below implying that the iterated kernels are symmetric measures, see Matthes et al.~\cite{MWM79}; sufficient conditions are given in Nehring and Zessin~\cite{NZ11}, Rafler~\cite{mR09} and Zessin~\cite{hZ09}.

\section{Characterization of Poisson \\and P\'olya Processes\label{sect:char}}

Led by the construction of Papangelou processes in Zessin~\cite{hZ09}, we fundamentally work under the following conditions on $\pi$ throughout the following sections:
\begin{enumerate}[leftmargin=5Em, labelsep=1.5Em]
  \item[$(\A_1)$] $\forall\mu\in\MpmX$: $\pi(\mu+\delta_y,\d x)\pi(\mu,\d y)=\pi(\mu+\delta_x,\d y)\pi(\mu,\d x)$,
  \item[$(\A_3)$] for $\mu\in\MpmX$, $B\in\Bbd$,
	\begin{equation*}
		0 < \sum_{m=0}^\infty \frac{\pi^{(m)}(\mu_{B^c},B^m)}{m!} < \infty,
	\end{equation*}
  \item[$(\A_7)$] $\pi$ is vagely continuous.
\end{enumerate}

The local integrability condition $(A_3)$ is important as soon as one is interested in the nonemptiness of $\wp(\pi)$. $\pi$ is vaguely continuous if $\pi$ is continuous when considered as a mapping $\pi:\MpmX\to\MX$ and both spaces are equipped with the vague topology. The vague continuity excludes kernels like $Z(\mu)\vo$, where $Z$ is a nonconstant function of the density of points, i.e. the avarage number of points per volume: One may construct a sequence of point measures starting from some fixed point measure $\mu_1$ with density of points $U(\mu_1)$, and add in each step one point such that the sequence converges vaguely to some point measure with a larger density of points; but since each $\mu_k$ differs from $\mu_1$ by a finite number of points, their densities agree.

What remains are postulates which relate the assumptions in Zessin~\cite{hZ09} to versions of Johnson's, Carnaps and B\"oge's ideas adapted to the point process setup. We start with a weak sufficiency postulate, consider then action of transformations of random measures and then turn to the learn-merge invariance principle.

\subsection{The weak sufficiency postulate}

The following condition is a special case of the assumptions in Zessin~\cite{hZ09} and the starting point of our reasoning:

\begin{enumerate}[leftmargin=5Em, labelsep=1.5Em]
  \item[$(\A_2)$] $\forall y\in X$  $\forall\mu\in \MpmX:$ $\pi(\mu+\delta_y,\d x) = \pi(\mu,\d x)$ on $\{y\}^c$. 
\end{enumerate}
$(\A_2)$ demands that adding a point at some location shall not be influenced by what happens outside that point. 
Condition $(\A_2)$ seems quite restrictive, but any density which is independent of $y$ needs to be constant due to the cocycle condition. We refer to the more general statement in Lemma~\ref{thm:intact:symm} below.

\prop{ \label{thm:cpp:char}
Under conditions $(\A_1),(\A_2)$ and $(\A_7)$ there exist a Radon measure $\vo\in\MX$ and a measurable function $c:X\times\N\to\R$ such that
\begin{equation}
  \pi(\mu,\,\cdot\,)=\vo+\sum_{x\in\supp\mu}c_x\bigl(\mu(x)\bigr)\delta_x
\end{equation}
In particular, $\vo = \pi(0,\,\cdot\,)$ and $c$ is the (signed) density $c_x(n)=\pi(n\delta_x,\{x\})-\pi\bigl((n-1)\delta_x,\{x\}\bigr)$.
}

The key to this theorem will be an induction argument based an $(\A_2)$ in connection with a comparison of the kernel $\pi$ evaluated at point configurations $\mu$ and $\nu$ which agree at almost all points. We postpone the proof and all following proofs to Section~\ref{sect:proofs}. 

The construction techniques developed in Nehring and Zessin~\cite{NZ11} and Zessin~\cite{hZ09} ensure that under the additional assumption of integrability, there exists a Papangelou process for a given kernel $\pi$.

\satz{ \label{thm}
Assume that in addition to the assumptions of Theorem~\ref{thm:cpp:char} the kernel $\pi$ also satisfies $(\A_3)$. Then there exists a unique Papangelou process $\Pp$ on $X$ with Papangelou kernel $\pi$. 
}

Johnson's original sufficiency postulate states that the probability of the observation of one of the finite number of categories shall only depend on the number of observations of that category (and the total number of observations so far). Zabell shows that one may allow this probability also to depend on the category itself. They restrict to the finite and the countably infinite case, respectively.

Suppose that for a kernel $\pi$, the intensity $\pi(\mu,\{x\})$ depends only on $x$ and $\mu(\{x\})$, i.e.
\begin{itemize}[leftmargin=5Em, labelsep=1.5Em]
  \item[$(\J)$]  $\forall x\in X:$ $\pi\bigl(\,\cdot\,,\{x\}\bigr)$ is $\sigma(\zeta_{\{x\}})$-measurable.
\end{itemize}

Alternatively, one may postulate that there exists a nonnegative function $F:X\times\N\to\R$ such that for all $ x\in X$ and $\mu\in\MpmX$, $\pi(\mu,\{x\})= F(x,\mu(\{x\}))$.

\satz{ \label{thm:cpp:john-equiv}
Assume that the Radon kernel $\pi$ satisfies the conditions $(\A_1)$ and $(\A_7)$. Then
\begin{enumerate}
	\item $(\A_2)$ implies $(\J)$,
	\item if $X$ is in addition discrete, then $(\J)$ implies $(A_2)$.
\end{enumerate}
}

\bem{
\begin{enumerate}
	\item There is no normalization depending on the number of realized points in contrast to the urn model context. 
	\item The reinforcement strength $c$ is local, i.e. is allowed to depend on the location and on the number of present points at that location but not on the whole configuration. In case of discrete $X$ this contrasts~\cite{sZ82,BM86} and is due to the fact that $\pi$ is not normalized. 
\end{enumerate}
}

\subsection{Transformations of Papangelou processes}

Let $Y$ be another c.s.m.s. and $T:\MX\to\M(Y)$ be a measurable mapping. Then the image measure of a random measure or point process $\Pp$ under $T$ is well defined and given by
\begin{equation*}
  T\Pp(\vi)=\Pp(\vi\circ T) = \int \vi(T\mu)\,\,\Pp(\d\mu)
\end{equation*}
for $\vi\in F\bigl(\M(Y)\bigr)$. $T\Pp$ is the law of $\Pp$ under $T$.

Of our particular interest is the following pointwise transformation: Let $G:X\to Y$ be a measurable mapping. $G$ is called \emph{proper} if
\begin{enumerate}[leftmargin=5Em, labelsep=1.5Em]
  \item[$(\pe)$] $\forall B\in\Bbd(Y):G^{-1}(B)\in\Bbd(X)$,
\end{enumerate}
i.e. pre-images of bounded sets are bounded themselves. The collection of these transformations we denote by $\Im(X,Y)$.

Let $G\mu$ be the image of $\mu\in\M(X)$ under $G$, i.e.
\begin{equation}
 G\mu(B) = \mu\bigl(G^{-1}(B)\bigr),\qquad B\in\B(Y).
\end{equation}
Then $\mu\to G\mu$ is a measurable mapping of $\M(X)$ into $\M(Y)$. 

\defn{State space transformation}{
Let $G:X\to Y$ be a measurable, proper mapping and $\Pp$ be a random measure on $X$. Then the distribution $G\Pp$ is called \emph{state space transformation} of $\Pp$.
}

Of which nature is property of being a Papangelou process under state space transformations? Let $\pi$ be a Papangelou kernel and $\Pp\in\wp(\pi)$ be a Papangelou process specified by $\pi$. The property of being a Papangelou process is preserved under a measurable, proper mapping $G:X\to Y$, if $\pi$ satisfies Dynkin's condition,
\begin{enumerate}[leftmargin=4.5Em, labelsep=1.5Em]
	\item[$(\D)$] $\forall\mu_1,\mu_2\!\in\!\MX\forall f\!\in\!\K(Y)\!:G\mu_1\!=\!G\mu_2 \Rightarrow \pi(\mu_1,\!f\circ G)\!=\!\pi(\mu_2,\!f\circ G)$,
\end{enumerate}
i.e. if the image of the two measures $\mu_1$ and $\mu_2$ under $G$ agree, then also the images of the kernel shall do. In this case, the kernel $\pi'$ from $\M(Y)$ to $Y$ given by
\begin{equation*}
  \pi'(\nu,B)=\pi(\mu,G^{-1}B),
\end{equation*}
where $\mu\in G^{-1}\nu$, is a well-defined Radon kernel. We observe that if a Papangelou kernel $\pi$ satisfies $(\D)$, then $\pi'$ is also a Papangelou kernel. This result is a version of Dynkin's state space transformation theorem for Markov processes~\cite{D} adapted to the Papangelou process context. For a related problems to Gibbs processes see Karrat and Zessin~\cite{K}.

\satzn{State space transformation}{ \label{thm:suff:sptt}
Let $G:X\to Y$ be a measurable, proper mapping of the c.s.m.s. $X$ onto the c.s.m.s. $Y$ and assume that $\pi(\mu,\,\cdot\,)$ is a Papangelou kernel on $X$ which satisfies $(\D)$. If $\Pp\in\wp(\pi)$, then $G\Pp\in\wp(\pi')$.
}
The state space transformation theorem has two immediate consequences on the set $\wp(\pi)$.
\korollar{
Under the condtions of Theorem~\ref{thm:suff:sptt},
\begin{equation*}
  G\wp(\pi)\subseteq\wp(\pi')
\end{equation*}
without equality in general.
}

\korollar{
Under the assumtions of the state space transformation theorem suppose that $Y=X$ and $G:X\to X$ leaves $\pi$ invariant, i.e. $G\pi = \pi$, then $\wp(\pi)$ remains invariant under $G$, i.e.
\begin{equation}
  G\wp(\pi) \subseteq \wp(\pi).
\end{equation}
In particular, if $\wp(\pi)=\{\Pp\}$ is a singleton, then $G\Pp = \Pp$.
}
\bem{
Poisson, P\'olya sum and P\'olya difference process given in example~\ref{ex:pap:kerne} are invariant under $G$ if and only if $G\vo=\vo$. This is a manifestation of the fact that they are characterized by their kernels. We remark also that in these three examples, the kernels $\pi$ even have the remarkable invariance property that they satisfy Dynkin's condition for every transformation $G\in \Im(X)=\cup_Y \Im(X,Y)$.
}

\subsection{Transformations and B\"oge's Learn-merge invariance}

In the terminology of B\"oge, an element $G\in\Im(X,Y)$ is called \emph{merger}, whereas the transformation $\Pp\mapsto G\Pp$ is the \emph{process of merging}, which means to coarsen observation categories. On the other hand, mapping a Papangelou process $\Pp$ to its Papangelou kernel $\pi_\Pp$ is a kind of \emph{process of learning} from the observation of a point configuration. The state space transformation theorem states that under Dynkins condition both processes commute, i.e. that the following diagrams commute for Papangelou processes $\Pp$:
\begin{equation}\label{diag:sptt:comp}
\begin{CD}
  \Pp		@>>>	\pi_\Pp\\
  @VVGV		        @VVGV\\
  G\Pp		@>>>	\pi'|G\pi_\Pp
\end{CD}\end{equation}
In this case, we say that the property of $\Pp$ being a Papangelou process is compatible with $G$.
\defn{$G$-compatibility}{
Let $G:X\to Y$ be a proper, measurable transformation. A Papangelou process $\Pp$ is \emph{$G$-compatible} if the diagramme~\eqref{diag:sptt:comp} commutes.
}
Indeed, if the Papangelou process $\Pp$ with Papangelou kernel $\pi$ is $G$-compatible, then $\pi$ satisfies Dynkins condition.
\prop{ \label{thm:suff:dynkin-equiv}
Let $G:X\to Y$ be a measurable, proper mapping of the c.s.m.s. $X$ onto the c.s.m.s. $Y$. If $\Pp$ is a Papangelou process with kernel $\pi$, then the following statements are equivalent:
\begin{enumerate}
  \item $\pi$ satisfies $(\D)$
  \item $\pi$ is $G$-compatible.
\end{enumerate}
}

The idea behind B\"oge's learn-merge invariance principle is that \emph{it should make no difference in inference, whether one combines (merges) some of the observation categories before or after sampling}~\cite{BM86}. For Papangelou processes, we study the implications of $G$-compatibility for a large class of proper mappings.

\begin{itemize}[leftmargin=5Em, labelsep=1.5Em]
 \item[(${\cal BC}$)]  for every continuous $G\in \Im(X)$, $\pi$ satisfies $(\D)$.
\end{itemize}

This strong condition implies the first version of Johnson's sufficiency postulate $(\J)$, i.e. $\pi(\mu,\{x\})$ depends only on $x$ and the multiplicity of $\mu$ at $x$, as we see in Theorem~\ref{thm:suff:CimplJ} below. Moreover, we show that $({\cal BC})$  is equivalent to a strong sufficiency postulate,
\begin{enumerate}[leftmargin=5Em, labelsep=1.5Em]
  \item [$(\J')$]  for every closed $B\in\Bbd(X)$, $\pi(\,\cdot\,,B)$ is $\sigma(\zeta_B)$-measurable.
\end{enumerate}

\satz{ \label{thm:suff:CimplJ}
For a Radon kernel $\pi$ satisfying $(\A_1)$, $({\cal BC})$ and $(\J')$ are equivalent. Moreover, both conditions imply $(\A_2)$ and thus $(J')$.
}

\bem{
A short example reveals the difference between $(\J)$ and $(\J')$: Let $X=\{-1,0,1\}$ and
\begin{equation*}
  \pi(\mu,\,\cdot\,)=\tfrac{1}{2}\left[1-\mu(\{-1\})\right]\delta_{-1}+\tfrac{1}{2}\delta_0+\tfrac{1}{2}\left[1+\mu(\{1\})\right]\delta_{1}.
\end{equation*}
Clearly, $\pi(\mu,\,\cdot\,)$ satisfies $(\J)$, but $\pi(\,\cdot\,,X)$ is not $\zeta_{X}$-measurable.
}

\satz{ \label{thm:suff:equiv}
Assume that $X$ contains at least 2 elements and let $\pi$ be a Radon kernel such that $(\A_1)$ and $(\A_7)$ are satisfied, then the following statements are equivalent:
\begin{enumerate}
  \item $\pi=\vo+c\mu$ for some $\vo\in\MX$ and constant $c$;
  \item $\pi$ satisfies $({\cal BC})$;
  \item $\pi$ satisfies $(\J')$.
\end{enumerate}
}

\bem{
Assuming the local integrability condition on $\pi$, which implies $c<1$, and in case $c<0$ that $\vo$ is a discrete measure with $\tfrac{\vo}{|c|}\in\Mpm(X)$, $(\J')$ as well as $({\cal BC})$ imply that $\wp(\pi)$ consists of a single element $\Pp$, which is either a Poisson, a P\'olya sum, or a P\'olya difference process.
}

\section{Interaction Included \label{sect:II}}

The conditions imposed on $\pi$ in Zessin~\cite{hZ09} did not exclude pair interactions as in Section~\ref{sect:char}. A generalization of the ideas that led to Theorem~\ref{thm:cpp:char} still allow to identify the basic structure of the kernel $\pi$ including a pair interaction. Firstly, we replace $(\A_2)$ by the weaker condition $(\A_2')$, and give additional conditions we need.
\begin{enumerate}[leftmargin=5Em, labelsep=1.5Em]
  \item[$(\A_2')$] there exists a measurable, nonnegative function $f$ on $X\times X$ such that $\pi(\mu+\delta_y,\d x) = f(x,y) \pi(\mu,\d x)$ on $\{y\}^c$ for each $y\in X$ and $\mu\in \MpmX$
  \item[$(\A_6)$] $f$ is stricly positive on the diagonal
\end{enumerate}

The strict positivity in $(\A_6)$ excludes a large class of interesting interactions as remarked at Example~\ref{ex:intact:boltzmann} below. We discuss the case of $f$ vanishing on the diagonal below. Condition $(\A_2')$ was originally assumed in Zessin~\cite{hZ09}. From the point of view of $(\A_2')$ we see the reason for the choice of a constant density in the last section.

\lemma{ \label{thm:intact:symm}
If the Radon kernel $\pi$ satisfies conditions $(\A_2')$, then $(\A_1)$ is equivalent to $f$ being a symmetric function.
}

\begin{proof}[Proof of Lemma~\ref{thm:intact:symm}]
By $(\A_1)$ and $(\A_2')$ we have
\begin{align*}
  1_{x\neq y}\pi(\mu+\delta_y,\d x)\pi(\mu,\d y) &=1_{x\neq y}\pi(\mu+\delta_x,\d y)\pi(\mu,\d x)\\
  1_{x\neq y}f(x,y)\pi(\mu,\d x)\pi(\mu,\d y) &=1_{x\neq y}f(y,x)\pi(\mu,\d y)\pi(\mu,\d x),
\end{align*}
hence
\begin{equation*}
  f(x,y)=f(y,x)\qquad\pi(\mu,\,\cdot\,)^2\text{-a.s. }(x,y).
\end{equation*}
\end{proof}

The induction argument used to prove Theorem~\ref{thm:cpp:char} suggests to define formally
\begin{equation*}
  V(x,\mu)=\prod_{z\in\supp\mu} f(x,z)^{\mu(\{z\})}.
\end{equation*}

\beispiel{\label{ex:intact:boltzmann}
Since $f$ is nonnegative, define $\phi(x,y)=-\log f(x,y)$ with $\phi(x,y)=+\infty$ if $f(x,y)=0$. Symmetry of $f$ implies symmetry of $\phi$, and formally 
\begin{equation*}
  V(x,\mu)=\exp\left(-\int \phi(x,z)\,\mu(\d z)\right).
\end{equation*}
Whenever this expression makes sense, the integral is the energy of a particle at $x$ given a configuration $\mu$, and $V$ is the corresponding Boltzmann factor.
}

Hence, assumption $(\A_6)$ excludes e.g. hard-core interactions and interactions via a repelling potential.

\satz{ \label{thm:intact:char}
Under conditions $(\A_1)$, $(\A_2')$, $(\A_6)$ and $(\A_7)$ there exist a Radon measure $\vo\in\MX$ and a measurable function $c:X\times\N\to\R$ such that
\begin{equation}
  \pi(\mu,\,\cdot\,)=V(x,\mu) \left(\vo+\sum_{x\in\supp\mu}c_x\bigl(\mu(\{x\})\bigr)\delta_x\right).
\end{equation}
In particular, $\vo = \pi(0,\,\cdot\,)$ and $c$ is the (signed) density
\begin{equation*}
	c_x(n)=f(x,x)^{-1}\pi(n\delta_x,\{x\})-\pi\bigl((n-1)\delta_x,\{x\}\bigr).
\end{equation*}
}

The first part is a classical interaction containing the Boltzmann-factor $V$. The second part is again a local reinforcement. If we know that $\pi$ is simple, then the reinforcement part vanishes and what remains is the Gibbs kernel.

\korollar{ \label{thm:intact:gibbs}
In addition to the assumptions of Theorem~\ref{thm:intact:char} assume that $\pi$ is simple. Then $\pi$ is the Gibbs kernel.
}

The condition of the positivity of $f$ in Theorem~\ref{thm:intact:char} means that $\phi$ needs to be finite even on the diagonal. At least a version of Theorem~\ref{thm:intact:char} remains true if $f$ vanishes on the diagonal.

\satz{ \label{thm:intact:pos}
Assume that the Radon kernel $\pi$ satisfies $(\A_1)$, $(\A_2')$ and $(\A_7)$ and the density $f$ vanishes on the diagonal, then there exist a Radon measure $\vo\in\MX$ and a measurable funtion $c:X\times\N\to\R$ such that
\begin{equation}
  \pi(\mu,\d z)=V(z,\mu)\vo(\d z)+\sum_{x\in\supp\mu} V(x,\mu|_{\{x\}^c}) c_x\bigl(\mu(\{x\})\bigr)\delta_x(\d z).
\end{equation}
}

Again the kernel is the sum of two parts, the first one is the Boltzmann kernel and the second one some reinforcement depending on the configuration some complicated way. However, if we know that $\pi$ is simple Corollary~\ref{thm:intact:gibbs} is still true under the assumption of $f$ vanishing on the diagonal.

\bem{
While in Section~\ref{sect:char} the additional assumption on the summability of the iterated kernels ensured existence and uniqueness of the corresponding point processes, the situation is different here.}
  
\section{Proofs \label{sect:proofs}}
\subsection{Proofs of Section~\ref{sect:char}}

We first prove Theorem~\ref{thm:cpp:char}, that the conditions $(\A_1)$, $(\A_2)$ and $(\A_7)$ imply that $\pi(\mu,\,\cdot\,)=\vo+\sum_{x\in\supp\mu} c_x\bigl(\mu(x)\bigr)\delta_x$.

\begin{proof}[Proof of Theorem~\ref{thm:cpp:char}]
From $(\A_1)$ and $(\A_2)$ we obtain the fundamental recursion
\begin{equation} \label{eq:fundrek}
  \pi(\mu+\delta_y,\,\cdot\,) = \pi(\mu,\,\cdot\,)+\bigl[\pi(\mu+\delta_y,\{y\})-\pi(\mu,\{y\})\bigr]\delta_y
\end{equation}
for $\mu\in\MpmX$ and $y\in X$. Setting $\vo=\pi(0,\,\cdot\,)$ as well as $c_x(1)=\pi(\delta_x,\{x\})-\vo(\{x\})$, equation~\eqref{eq:fundrek} specializes to
\begin{equation}
  \pi(\delta_y,\,\cdot\,) = \vo+c_y(1)\delta_y,
\end{equation}
thus
\begin{equation} \label{eq:claim}
  \pi(\mu,\,\cdot\,) = \vo+\sum_{x\in\supp\mu} c_x\bigl(\mu(x)\bigr)\delta_x
\end{equation}
holds for all point configurations $\mu\in\Mpm_{1}(X)$. Next, we show inductively that equation~\eqref{eq:claim} also holds for all $\mu\in\Mpm_n(X)$ for all $n\in\N$.

Assume that~\eqref{eq:claim} holds for all $\mu\in\Mpm_n(X)$ for some $n\geq 1$. Together with equation~\eqref{eq:fundrek} this implies
\begin{equation*}
  \pi(\mu+\delta_y,\,\cdot\,)=\vo+\sum_{x\in\supp\mu} c_x\bigl(\mu(x)\bigr)\delta_x+ F(\,\cdot\,,\mu)\delta_y,
\end{equation*}
where
\begin{equation*}
  F(x,\mu)=\pi(\mu+\delta_x,\{x\})-\vo(\{x\})-\sum_{x\in\supp\mu} c_x\bigl(\mu(x)\bigr)\delta_x.
\end{equation*}

If $\mu=n\delta_y$, then we are done by setting $c_y(n+1)=c_y(n)+F(y,\mu)$ and by observing that $F(y,\mu)$ depends on $\mu$ only via $\mu(y)$. Otherwise, choose $x\in\supp\mu\cap\{y\}^c$, then also $\nu=\mu-\delta_x+\delta_y\in\Mpm_n(X)$, $\mu$ and $\nu$ agree on $\{x,y\}^c$, and 
\begin{align*}
  \pi(\mu+\delta_y,\,\cdot\,) &= \vo +\sum_{\substack{z\in\supp\mu\\ z\neq x,y}} c_z\bigl(\mu(z)\bigr)\delta_z +c_x\bigl(\mu(x)\bigr)\delta_x +c_y\bigl(\mu(y)\bigr)\delta_y +F(\,\cdot\,,\mu)\delta_y\\
  \pi(\nu+\delta_x,\,\cdot\,) &= \vo +\sum_{\substack{z\in\supp\mu\\ z\neq x,y}} c_z\bigl(\nu(z)\bigr)\delta_z + c_x\bigl(\mu(x)-1\bigr)\delta_x +c_y\bigl(\mu(y)+1\bigr)\delta_y\\
  	 &\qquad+F(\,\cdot\,,\nu)\delta_x.
\end{align*}

Thus, by comparing coefficients, we get
\begin{equation*}
	c_y\bigl(\mu(y)+1\bigr)= c_y\bigl(\mu(y)\bigr)+F(y,\mu),
\end{equation*}
and therefore~\eqref{eq:claim} holds on the set of all finite point measures.

Finally, note that the set of finite point measures is a vaguely dense subset of $\MpmX$, and therefore, by the vague continuity of the kernel, equation~\eqref{eq:claim} extends to $\MpmX$.
\end{proof}

\begin{proof}[Proof of Theorem~\ref{thm:cpp:john-equiv}]
By Theorem~\ref{thm:cpp:char}, $(\A_2)$ clearly implies $(\J)$. Next assume that there exists a measurable function $g$ such that $\pi(\mu,\{x\})=g\bigl(\mu(x),x\bigr)$. If $X$ is discrete and $B\in\B$ is chosen such that $y\notin B$, then
\begin{equation*}
	\pi(\mu+\delta_y,B)=\sum_{x\in B} \pi(\mu+\delta_y,\{x\})=\sum_{x\in B} g\bigl(\mu(x),x\bigr)=\pi(\mu,B).
\end{equation*}
Thus $\pi(\mu+\delta_y,\,\cdot\,)\ll\pi(\mu,\,\cdot\,)$ on $\{y\}^c$ with density 1.
\end{proof}

Next we show that under Dynkin's condition, a mapped Papangelou process is again a Papangelou process
\begin{proof}[Proof of Theorem~\ref{thm:suff:sptt}]
Firstly, note that by Dynkin's condition and surjectivity of $T$,  $\pi'$ is a well-defined Radon kernel on $Y$.

The following calculation shows that if $\Pp\in\wp(\pi)$, then $G\Pp$ is a Papangelou process with Papangelou kernel $\pi'$. Let $h\in F(Y\times \M(Y))$, then
\begin{align*}
  C_{G\Pp}(h) &= \iint_{Y\times\M(Y)} h(y,\nu)\,\nu(\d y)G\Pp(\d\nu)\\
	&= \iint_{X\times\M(X)} h(Gx,G\mu)\,\mu(\d x)\Pp(\d\mu)\\
	&= \iint_{X\times\M(X)} h\bigl(Gx,G(\mu+\delta_x)\bigr)\,\pi(\mu,\d x)\Pp(\d\mu).
\intertext{Now observe that $G(\mu+\delta_x)=G\mu + \delta_{Gx}$, and therefore}
	&= \iint_{X\times\M(X)} h(Gx,G\mu+\delta_{Gx})\,\pi(\mu,\d x)\Pp(\d\mu).
\intertext{Finally, Dynkin's condition allows to transform the integral back to $Y$,}
	&= \iint_{Y\times\M(Y)} h(y,\nu+\delta_{y})\,\pi'(\nu,\d y) G\Pp(\d\nu),
\end{align*}
hence $G\Pp\in\wp(\pi')$.
\end{proof}

\begin{proof}[Proof of Proposition~\ref{thm:suff:dynkin-equiv}]
Suppose that the diagramme~\eqref{diag:sptt:comp} commutes, which means that for any $h\in F(Y\times\M(Y))$,
\begin{equation*}
  \iint h(y,G\mu+\delta_y)\pi'(G\mu,\d y)\Pp(\d\mu)=\iint h(Gx,G\mu+\delta_{Gx})\pi(\mu,\d x)\Pp(\d\mu).
\end{equation*}
Choose $h=\exp(-f)\otimes\exp(-\zeta_g)$ for given $f,g:Y\to\R$ with bounded support, then
\begin{equation*}
  \pi'(G\mu,\exp(-f-g))=\pi(\mu,\exp(-f\circ G-g\circ G))
\end{equation*}
for $\Pp$-a.e. $\mu$ and since the lhs depends only on $G\mu$, $\pi(\,\cdot\,,\exp(-f\circ G-g\circ G))$ is a.s. constant on $\{G=\nu\}$ and $\pi$ satisfies Dynkin's condition for $G$.
\end{proof}

Next we show that $({\cal BC})$ and $(\J')$ are equivalent. Recall that the underlying space $X$ is equipped with a metric $d$.

\begin{proof}[Proof of Theorem~\ref{thm:suff:CimplJ}]
Firstly assume $({\cal BC})$. Let $B\in\Bbd(X)$ be closed. Define $G:X\to[0,\infty)$, $x\mapsto d(x,B)$. Since $B$ is bounded, $G$ is proper, it is clearly continuous (and thereby measurable) and since $B$ is closed, $x\in B $ holds if and only if $G(x)=0$. 

Next choose $\mu_1,\mu_2$ such that $G\mu_1=G\mu_2$. Then by assumption, $\pi(\mu_1,B)=\pi(\mu_2,B)$. On the other hand, $G\mu_1=G\mu_2$ implies $\mu_1(B)=\mu_2(B)$, and therefore, $\pi(\,\cdot\,,B)$ is measurable with respect to $\sigma(\zeta_B)$, which is $(\J')$. Specializing to singletons $B$ then imples $(\J)$.

Now assume that $(\J')$ holds and follow the standard construction of integrals. Let $G:X\to Y$ be a proper mapping and choose $\mu_1$ and $\mu_2$ such that $G\mu_1=G\mu_2$. We have to show that
\begin{equation*}
  \pi(\mu_1,g\circ G)=\pi(\mu_2,g\circ G)
\end{equation*}
for all measurable $g:Y\to\R_+$. Choose $A\in\Bbd(Y)$, then since $G\mu_1=G\mu_2$, $G\mu_1(A)=G\mu_2(A)$. Moreover,
\begin{equation*}
  \pi(\mu_1,1_A\circ G)=\pi(\mu_1,G^{-1}(A))
\end{equation*}
is $\sigma(\zeta_{G^{-1}(A)})$-measurable by assumption, hence
\begin{equation*}
  \pi(\mu_1,G^{-1}(A))=\pi(\mu_2,G^{-1}(A))=\pi(\mu_2,1_A\circ G).
\end{equation*}
By linearity and monotone convergence, we obtain the result for any $g$, hence $({\cal BC})$ holds.

Finally, to show that $(\J')$ implies $(\A_2)$, let $y\in X$ and $B\in\F_{\{y\}^c}$. Then
\begin{equation*}
	\pi(\mu+\delta_y,B)=g\bigl(B,\mu(B)\bigr)
\end{equation*}
since $y\notin B$. But the r.h.s. equals $\pi(\mu,B)$.
\end{proof}

We know that $(J')$ implies that the kernel $\pi$ is composed of a fixed measure $\vo$ and a local reinforcement. Hence, of interest is the structure of the reinforcement.

\begin{proof}[Proof of Theorem~\ref{thm:suff:equiv}]
Clearly a linear structure of $\pi$ implies $(\J')$ and $({\cal BC})$.

To show that the strong sufficiency postulate implies the linear structure, fix $B\in\B_0$ such that $B$ contains at least 2 elements. By choosing configurations which are concentrated on a single element in $B$, we see that for each $n$ the function $x\mapsto c_x(n)$ is constant. Moreover, by induction over the number of points of configurations, we get that $c(n)=n\cdot c(1)$.
\end{proof}

\subsection{Proofs of Section~\ref{sect:II}}

The main problem here is to take care of the density $f$. Fortunately, the arguments carry over mostly.

\begin{proof}[Proof of Theorem~\ref{thm:intact:char}]
In analogy to~\eqref{eq:fundrek}, we obtain from $(\A_1)$ and $(\A_2')$ the fundamental recursion
\begin{equation} \label{eq:G:fundrek}
  \pi(\mu+\delta_y,\,\cdot\,) = f(\,\cdot\,,y)\pi(\mu,\,\cdot\,)+\bigl[\pi(\mu+\delta_y,\{y\})-f(y,y)\pi(\mu,\{y\})\bigr]\delta_y
\end{equation}
for $\mu\in\M$ and $y\in X$. Define $\vo=\pi(0,\,\cdot\,)$ and $c_x(1) =\tfrac{1}{f(x,x)}\pi(\delta_x,\{x\}) -\pi(0,\{x\})$, then
\begin{equation} \label{eq:proofs:f-rausziehen}
  \pi(\delta_y,\,\cdot\,) = f(\,\cdot\,,y)\bigl[\vo+c_y\delta_y\bigr].
\end{equation}

Thus, if we denote by $\pi_0$ the kernel identified in Theorem~\ref{thm:cpp:char}, 
\begin{equation}\label{eq:G:claim}
  \pi(\mu,\,\cdot\,) = V(\,\cdot\,,\mu)\pi_0(\mu,\,\cdot\,) 
\end{equation}
for all $\mu\in\Mpm_1(X)$. Assuming that~\eqref{eq:G:claim} holds for all $\mu\in\M_n$ for some $n\in\N$, by using~\eqref{eq:G:fundrek} we get
\begin{equation*}
  \pi(\mu+\delta_y,\d z) = V(z,\mu+\delta_y)\pi_0(\mu,\d z)+ F(z,y,\mu)\delta_y(\d z),
\end{equation*}
where
\begin{equation*}
  F(z,y,\mu)=\pi(\mu+\delta_y,\{z\})-f(z,z)V(z,\mu)\pi_0\bigl(\mu,\{z\}\bigr).
\end{equation*}

From this point, the main arguments agree with those of the proof of Theorem~\ref{thm:cpp:char} with the added Boltzmann factor $V$. We do not repeat them here.
\end{proof}

What happens if $f$ vanishes on the diagonal? Firstly, equation~\eqref{eq:proofs:f-rausziehen} takes a different form and already the choice choice of the function $c$ has to be adapted. We give comments on the changes due to the zeros. 

\begin{proof}
In case of $f$ vanishing on the diagonal, the fundamental recursion simplifies to
\begin{equation}
  \pi(\mu+\delta_y,\,\cdot\,) = f(\,\cdot\,,y)\pi(\mu,\,\cdot\,)+\pi(\mu+\delta_y,\{y\})\delta_y.
\end{equation}
$\vo$ is defined analogously and $c_y(n)=\pi(n\delta_y,\{y\})$. Next, the induction step turns into
\begin{align*}
	\begin{multlined}
		\pi(\mu+\delta_y,\d z) = V(z,\mu+\delta_y)\vo(\d z)\\ +\sum_{x\in\supp\mu} f(z,y)V(x,\mu|_{\{x\}^c})c_x\bigl(\mu(\{x\})\bigr)\delta_x(\d z)
	 +\pi(\mu+\delta_y,\{y\})\delta_y(\d z).\end{multlined}
\end{align*}
Observe that we may replace $f(z,y)$ by $f(x,y)$ and moreover
\begin{equation*}
	f(x,y)V(x,\mu|_{\{x\}^c}) =\begin{cases}
		V(x,\mu|_{\{x\}^c}+\delta_y) &\text{ if } x\neq y\\
		0 &\text{ if } x=y.
	\end{cases}
\end{equation*}
If $\mu$ is concentrated on $\{y\}$, we are done; otherwise, by replacing $\mu+\delta_y=\nu+\delta_x$ we get
\begin{equation*}
	\pi(\mu+\delta_y,\{y\}) = f(y,x)V(y,\nu|_{\{y\}^c})c_y\bigl(\nu(\{y\})\bigr).
\end{equation*}
Finally, remark that $\nu|_{\{y\}^c}+\delta_x=\mu|_{\{y\}^c}$ and $\nu(\{y\})=\mu(\{y\})+1$, and this last observation yields the claim.
\end{proof}

\subsection*{Acknowledgements}
The authors thank the referee for carefully reading the manuscript and giving valuable comments.




\begin{thebibliography}{99}

\bibitem{BZ}
Bach, A., Zessin, H.: The particle structure of the quantum mechanical Bose and Fermi gas. Unpublished manuscript (2010).


\bibitem{BM86}
B\"oge, W., M\"ocks, J.: Learn-merge invariance of priors: A characterization of the Dirichlet distributions and processes. J. Multivar. Anal. 18, 83 - 92 (1986).

\bibitem{CS58}
Carnap, R., Stegm\"uller, W.: Induktive Logik und Wahrscheinlichkeit. Springer, Wien (1958).

\bibitem{D}
Dynkin, E.B., Markov Processes: vol.1 , Springer, Berlin (1965).

\bibitem{Ess}
Essler, W.K.: Induktive Logik, Alber-Verlag, Freiburg und M\"unchen (1970).

\bibitem{G76}
Georgii, H.-O.: Canonical and grand canonical Gibbs states for continuum systems, Commun. Math. Phys. 48, 31 - 51 (1976).

\bibitem{G}
Good, I.J.: The Estimation of Probabilities. MIT Press Cambridge (1967).

\bibitem{H}
Humburg, J.: Die Problematik apriorischer Wahrscheinlichkeiten im System der induktiven Logik von Rudolf Carnap. Arch. Math. Logik Grundlag 14 , 135 - 147 (1971).

\bibitem{J}
Johnson, W.E., Braithwaite, R.B.: Probability: The deductive and inductive problems (appendix). Mind 41 , 421 - 423 (1932).

\bibitem{Ka}
Kallenberg, O.: Random Measures, Akademie-Verlag and Academic Press, Berlin (1983).

\bibitem{K}
Karrat, M., Zessin, H.: Transformations of Gibbs states, CRAS 332, S\'erie I, 453 - 458 (2001).

\bibitem{MKM78}
Kerstan, J., Matthes, K. and Mecke, J.:
\newblock {Infinitely Divisible Point Processes}.
\newblock Wiley, Chichester (1978).


\bibitem{kK72}
Krickeberg, K.: The Cox Process, Symp. Math. 9, Calcolo Probab., Teor. Turbolenza 9, 151 - 167 (1972).


\bibitem{LMW88}
Liemant, A., Matthes, K., Wakolbinger, A.: Equilibrium distributions of Branching Processes, Akademie-Verlag, Berlin (1988).

\bibitem{MWM79}
Matthes, K., Warmuth, W. and Mecke, J.:
\newblock Bemerkungen zu einer Arbeit von Nguyen Xuan Xanh and Hans Zessin.
\newblock {\em Math. Nachr.} 88, 117--127 (1979).


\bibitem{NZ11}
Nehring, B., Zessin, H.: The Papangelou process. A concept for Gibbs, Fermi and Bose processes, J. Contemp. Math. Anal. 46, 326 - 337 (2011).

\bibitem{NZ79}
Nguyen, X.~X. and Zessin, H.:
\newblock Integral and differential characterisation of the Gibbs process.
\newblock {\em Math. Nachr.} 88, 105--115 (1976/77). 

\bibitem{P}
Papangelou, F.: Point processes on spaces of flats and other homogeneous spaces, Math. Proc. Camb. Philos. Soc. 80, 297 - 314 (1976).

\bibitem{P1}
Papangelou, F.: The conditional intensity of general point processes and an application to line processes. Z. Wahrscheinlichkeitstheorie Verw. Geb. 28, 207 - 226 (1974).

\bibitem{mR09}
Rafler, M.: Gaussian Loop- and P\'olya Processes. A Point Process Approach, Universit\"atsverlag Potsdam (2009).

\bibitem{R}
Rafler, M.: A Cox process representation for the P\'olya sum process, J. Contemp. Math. Anal. 46, 338 - 345 (2011).

\bibitem{sZ82}
Zabell, S.:
\newblock W. E. Johnson's "sufficientness" postulate.
\newblock {\em Ann. Stat.} 10(4), 1091--1099 (1982). 

\bibitem{hZ09}
Zessin, H., Der Papangelou Proze\ss{} , J. Contemp. Math. Anal. 44, 61 - 72 (2009).

\end{thebibliography}
\end{document}